\documentclass{llncs}
\usepackage{epsfig}

\newcommand{\A}{\mathcal{A}}

\newcommand{\R}{\mathcal{R}}

\newcommand{\T}{\mathcal{T}}

\newcommand{\ob}{\langle}
\newcommand{\cb}{\rangle}

\newcommand{\mn}{\hbox{min}}
\newcommand{\mx}{\hbox{max}}

\begin{document}

\pagestyle{headings}

\title{
New Developments in Interval Arithmetic and Their Implications for
Floating-Point Standardization\thanks{
Technical report DCS-273-IR, Department of Computer Science,
University of Victoria, Victoria, BC, Canada.}
}

\author{M.H. van Emden}

\institute{Department of Computer Science, \\University of Victoria,
Victoria, Canada \\
\email{vanemden@cs.uvic.ca}\\
\texttt{http://www.cs.uvic.ca/\homedir vanemden/}
}

\maketitle

\begin{abstract}
We consider the prospect of a processor that can perform interval
arithmetic at the same speed as conventional floating-point
arithmetic.  This makes it possible for all arithmetic to be performed
with the superior security of interval methods without any penalty in
speed. In such a situation the IEEE floating-point standard needs to be
compared with a version of floating-point arithmetic that is ideal for
the purpose of interval arithmetic.  Such a comparison requires a
succinct and complete exposition of interval arithmetic according to
its recent developments. We present such an exposition in this paper.
We conclude that the directed roundings toward the infinities and the
definition of division by the signed zeros are valuable features of the
standard.  Because the operations of interval arithmetic are
always defined, exceptions do not arise. As a result neither Nans nor
exceptions are needed. Of the status flags, only the inexact flag may
be useful. Denormalized numbers seem to have no use for interval
arithmetic; in the use of interval constraints, they are a handicap.

\end{abstract}


\section{Introduction}
Continuing advances in process technology have caused a tremendous
increase in the number of transistors available to the designer of a
processor chip. As a result, multiple parallel floating-point units
become feasible. The time will soon come when interval arithmetic can
be done as fast as conventional arithmetic.

However, to properly utilize the newly available number of transistors,
chip designers need to spend ever more time iterating through cycles of
synthesis, place-and-route, and physical verification than current
design methodology allows. This design bottleneck makes it desirable
to simplify floating-point units. Recent developments in the
theory of interval arithmetic suggest possibilities for simplification.
As far as interval arithmetic is concerned, certain parts the 1985
standard are essential, whereas other parts are superfluous or even a
liability.

In this paper we present a consolidated, self-contained account of the new
developments in interval arithmetic not available elsewhere. In the
conclusions we give a sketch of what would be an ideal standard from
the point of view of interval arithmetic when arranged in this way.

\section{Why interval arithmetic?}

Conventional, non-interval, numerical analysis is marvelously cheap and
it works most of the time.  This was exactly what was needed in the
1950s when computers needed a demonstration of feasibility. A lot has
changed since that time.  Numerical computation no longer needs to be
cheap. It has become more important that it always works. As a result
interval arithmetic is becoming an increasingly compelling
alternative.  For example, civil, mechanical, and chemical engineers are
liable for damage due to unsound design. So far, they have been able to
get away with the use conventional numerical analysis, appealing to
what appears to be best practice. As interval methods mature, it is
becoming harder to ignore them when defining ``best practice''.


Recent developments, which we call ``modern interval arithmetic''
provide a practical and mathematically compelling basis. In
conjunction with this it has become clear that some aspects of IEEE
standard 754 are not needed or are detrimental, whereas other aspects
are marvelously suited to interval arithmetic. If these latter are
preserved in future development of the standard, then interval
arithmetic can help bridge the design gap and lead to the situation
where all arithmetic can be faster due to interval methods.

\section{A theory of approximation}

Conventional numerical analysis approximates each real by a single
floating-point number. It approximates each of the elementary
operations on reals by the floating-point operation that has the same
name, but not the same effect. Let us call this approach ``point
approximation''. It has been amply documented that this approach,
though often satisfactory, can lead to catastrophical errors. This can
happen because it is not known to what degree a floating-point variable
approximates its real-valued counterpart.

Interval arithmetic is based on a theory of approximation, {\em
set approximation}, that ensures that for every real-valued variable
$x$ in the mathematical model, there is a machine-representable set $X$
of reals that contains $x$. Such arithmetic is exact in the
sense that $x \in X$ is and remains a true statement in the sense of
mathematics. It is of course not exact in the sense that $X$ typically
contains many reals.

Conversely, in case of numerical difficulties, it will turn out that
continued iteration does not reduce the size of $X$, in which case we
have a notification that there is a problem with the algorithm or with
the model. Because of this property, we call this {\em manifest}
approximation: there is always a known lower bound to the quality of
the approximation.

In this way the operations can be interpreted as
inference rules of logic; for example, $x \in X$ and $y \in Y$ imply
that $x+y \in Z$, where $Z$ is computed from $X$ and $Y$.

\section{Approximation structures}

We address the situation where we need to solve a mathematical model in
which a variable takes on values that are not representable in a
computer, but where it is possible to so represent {\em sets} of
values. We then approximate a variable by a representable set that
contains all the values that are possible according to the model.
Models with real-valued variables are but one example of such a
situation.

As the theory of set approximation applies to sets in general, we
first present it this way.

\begin{definition}
\label{approxStruct}
Let $\T$ be the type of a variable $x$, that is, the set of the
possible values for $x$.  A finite set $\A$ of subsets of $\T$ that contains
$\T$ and that is closed under intersection is called an {\em
approximation structure in $\T$}.
\end{definition}

In a typical practical application of this theory, $\A$ is
a set of computer-representable subsets of $\T$.

\begin{theorem}
\label{approxTh}
If $\A$ is an approximation structure for $\T$, then
for every $S \subset \T$ there exists a unique least
(in the sense of the set-inclusion partial order) element
$S'$ of $\A$ such that $S \subset S'$.
\end{theorem}

\begin{definition}
\label{phi}
For every $S \subset \T$, $\phi(S)$ is the unique least element of $\A$
that exists according to theorem~\ref{approxTh}.
\end{definition}
We regard $\phi(S)$ as the approximation of $S$.
As $S$ can be a singleton set, this theory
provides approximations both of elements and subsets of $\T$.

\section{An approximation structure for the reals}

We seek a set of subsets of the reals that can serve as approximation
structure. A first step, not yet computer-representable, is that of
the closed, connected sets of reals.

\begin{theorem}
\label{cloConn}
The closed, connected subsets of $R$ are an approximation
structure for $R$.
\end{theorem}

{\em Proof}
According to a well-known result in topology,
all closed connected subsets of $R$ have one of the following
four forms:
$\{x \in \R \mid a \leq x \leq b \}$,
$\{x \in \R \mid x \leq b \}$,
$\{x \in \R \mid a \leq x \}$, and
$\R$
where $a$ and $b$ are reals.
Here we do not exclude $a > b$, because the empty set is also
included among the closed connected sets.
Clearly the conditions for an approximation structure are satisfied:
$\R$ is included and the intersection of any two sets of this form is
again closed and connected.
This completes the proof.

The significance of the closed connected sets of reals as approximation
structure is that they can be represented as a pair of extended reals.

\begin{definition}
\label{extReals}
The {\em extended reals} are the set obtained by adding to the reals
the two infinities. As with the reals, the extended reals are totally
ordered. When two extended reals are finite, then they are ordered
within the extended reals as they are in the reals. Furthermore,
$-\infty$ is less than any real and $+\infty$ is greater than any
real.
\end{definition}

Theorem~\ref{cloConn} together with definition~\ref{extReals}
suggest the following notation for the four forms of the closed
connected sets of reals. In this notation we do not include the empty
interval. The reason is that if it is found in an interval computation
that an interval is empty, any further operations involving the
corresponding variable yield the same result, so that the computation
can be halted.

\begin{definition} \label{brackets}
Let $a$ and $b$ be reals such that $a \leq b$.
\begin{eqnarray*}
\ob a,b \cb & \stackrel{\rm def}{=} &
     \{x \in \R \mid a \leq x \leq b \} \\
\ob -\infty,b \cb & \stackrel{\rm def}{=} &
     \{x \in \R \mid x \leq b \}        \\
\ob a,+ \infty \cb & \stackrel{\rm def}{=} &
\{x \in \R \mid a \leq x \}             \\
\ob -\infty, +\infty \cb & \stackrel{\rm def}{=} & \R
\end{eqnarray*}
\end{definition}

Note that each of these pairs denote sets of reals, even though in
their notation the infinities are used. These are not reals.

\section{Floating-point intervals: a finite approximation structure for
the reals}

Let $F$ be a finite set of reals.

\begin{theorem}
\label{finApprox}
The sets of the form
$\emptyset$,
$\ob -\infty,b \cb$,
$\ob a,b \cb$,
$\ob a,+\infty \cb$, and
$\ob -\infty,+\infty \cb$
are an approximation structure when $a$ and
$b$ are restricted to elements of $F$ such that $a \leq b$.
\end{theorem}

\begin{definition}
\label{flptOrReal}
The {\em real intervals} are sets of the form described in
theorem~\ref{cloConn}.
The {\em floating-point intervals} are sets of the form described in
theorem~\ref{finApprox}, where $a$ and $b$ are finite IEEE754 floating
point numbers (according to a choice of format:  single-length,
double-length, extended) such that $a \not= -0$, $b \not= +0$, and $a
\leq b$.
\end{definition}

From definitions \ref{brackets}
and \ref{flptOrReal} and theorem \ref{finApprox} we conclude:

\begin{itemize}
\item
The restriction on the sign of zero bounds
in definition~\ref{flptOrReal} is there to make
the notation unambiguous. We will see that disambiguating the notation
in this way has an advantage for interval division.

\item
$\{0\}$ is written as $\ob+0,-0\cb$.

\item
When $\ob a,b \cb$ is a floating-point interval, then $a \not= +\infty$
and $b \not= -\infty$.
\end{itemize}

Let us take care to distinguish ``real intervals'' from
``floating-point intervals''. Both are sets of reals. The latter are a
subset of the former.

From now on we assume
the floating-point intervals
as approximation structure when we rely on
the fact that for any set $S$ of reals there is a unique least
floating-point interval $\phi(S)$ containing it.

\begin{definition}
For any real $x$, $x^-$ ($x^+$) is the left (right) bound
of $\phi(\{x\})$.
\end{definition}

This operation is implemented by performing a floating-point
operation that yields $x$ in rounding mode toward $-\infty$
($+\infty$).

\section{Interval Arithmetic}
Much of the standard is concerned with defining, signaling, and
trapping exceptions caused by overflow, underflow, and undefined
operations. What distinguishes modern interval arithmetic from the old
is that {\em no exceptions occur}.  As we will see, no operation can
result in Nan.  Every operation is defined on all operands. Moreover,
it is defined in such a way that the floating-point endpoints bound the
set of the real numbers that are the possible values of the associated
variable in the mathematical model.

This property is based on the use of {\em set extensions} of the
arithmetical operations. It is helped by the use of {\em relational
definition} rather than functional ones of these operations. We
discuss these in turn.

\paragraph{Set extensions of functions}
Whenever a function $f$ is defined on a set $S$ and has values in a set
$T$, there exists the {\em canonical set extension} $\widehat f$, which
is a function defined on the subsets of $S$ and has as values subsets
of $T$ according to $\widehat f(X) = \{ f(x) \mid x \in X\}$ for any $X
\subset S$. This definition is of interest because it also carries over
to partial functions and to multivalued functions.

Though $X$ may be an approximation of $x$, $\widehat f(X)$ may not be an
element of an approximation structure of $T$, so is not necessarily
an approximation of $f(x)$. But $\phi(\widehat f(X))$ does approximate $f(x)$.
Thus $\phi$ induces a transformation among functions.
It changes $f$ to
the function that maps $x$ to $\phi(\widehat f(\{x\}))$.

The {\em inverse canonical set extension} of $f$ is defined as
$f^{-1}(Y) = \{ x \mid f(x) \in Y \}$. This definition is of interest
because such an inverse is defined even when $f$ itself has no
inverse.

By using the canonical set extensions of a function, one ensures that
undefined cases never arise. By considering instead of the
arithmetical operations on the reals their canonical set extensions
to suitably selected sets of reals (namely, floating-point intervals),
undefined cases are eliminated.

An example of a set extension for arithmetical operations is
$X+Y = \{x + y \mid x \in X \wedge y \in Y\}$. Though $X$ and $Y$ may be
floating-point intervals, that is typically not the case for
$\{x + y \mid x \in X \wedge y \in Y\}$. So
to ensure that addition is closed in the set of floating-point
intervals, we need to apply $\phi$, as shown below in the formulas for
interval operations that go back to R.E. Moore \cite{moore66}.

\begin{eqnarray}
X+Y & = &
 \phi(\{x + y \mid x \in X \wedge y \in Y\}) \nonumber\\
X-Y & = &
 \phi(\{x - y \mid x \in X \wedge y \in Y\}) \nonumber\\
X*Y & = &
 \phi(\{x * y \mid x \in X \wedge y \in Y\}) \nonumber\\
X/Y & = &
 \phi(\{x / y \mid x \in X \wedge y \in Y\}) \nonumber
\end{eqnarray}        

Regarded as a set extension, the above definition of $X/Y$ is correct
and unambiguous: set extensions are defined just as well for partial
functions, functions that are not everywhere defined.  Yet many authors
have subjected it to the condition $0 \not\in Y$, making it useless in
practice.  Others have taken a less restrictive stance by changing the
definition to:
$$
X/Y = 
\phi(\{x / y \mid x \in X \wedge y \in Y \wedge y \not= 0\}).
$$

\paragraph{Relational definitions}

Ratz \cite{ratz96} has avoided such difficulties by using a relational form
of the above definitions. Although not necessary, this relational form
also makes it possible to define both addition and subtraction with the
same ternary relation $x+y=z$. This leads to an attractive
uniformity in the definition of the interval arithmetic operations.

\begin{definition}
\label{intArithm}
Let $X$ and $Y$ be non-empty floating-point intervals. Then interval
addition, subtraction, multiplication, and division are defined as
follows.
\begin{eqnarray}
X+Y & \stackrel{\rm def}{=} &
 \phi(\{z \mid \exists x \in X \wedge
                \exists y \in Y.\; x+y=z\}) \nonumber\\
X-Y & \stackrel{\rm def}{=} &
  \phi(\{z \mid \exists x \in X \wedge
                 \exists y \in Y.\; z+y=x\}) \nonumber\\
X*Y & \stackrel{\rm def}{=} &
  \phi(\{z \mid \exists x \in X \wedge
                 \exists y \in Y.\; x*y=z\}) \nonumber\\
X \oslash Y & \stackrel{\rm def}{=} &
  \phi(\{z \mid \exists x \in X \wedge
                 \exists y \in Y.\; z*y=x\}) \nonumber
\end{eqnarray}
\end{definition}

We use the symbol $\oslash$ in $X \oslash Y$ here for interval division
rather than the $X/Y$ defined earlier. There is only a difference
between the two definitions when $\ob 0,0 \cb$ occurs as an operand.
For details, see \cite{hckvnmdn01}.  The difference is immaterial, as
intuition fails in these cases, anyway.

The operations thus defined form an interval arithmetic that is {\em
sound} in the sense that the resulting sets contain all the real
values they should contain according the set extension definition.
They are {\em closed} in the sense that they are defined for all
interval arguments and yield only interval results. Such an interval
arithmetic never yields an exception.

It remains to show that these definitions can be efficiently computed
by IEEE standard floating-point arithmetic while avoiding the undefined
floating-point operations
$\infty - \infty$, $\pm \infty / \pm \infty$, $0 * \pm
\infty$, and $0/0$.  This we do in the next sections.

\subsection{The algorithm for interval addition and subtraction}

\begin{theorem}
If $X=\ob a,b\cb $ and $Y=\ob c,d\cb $ are non-empty
floating-point intervals, then
$X+Y$ and $X-Y$ according to definition~\ref{intArithm} are equal to
$\ob (a+c)^-, (b+d)^+ \cb$ and
$\ob (a-d)^-, (b-c)^+\cb$, respectively. 
\end{theorem}

See \cite{hckvnmdn01}.
The interesting part of the proof takes into account that
adding $a$ and $c$ is undefined if they are infinities with opposite
signs. As, according to definition~\ref{brackets},
$a$ and $c$ are not $+\infty$, this cannot
happen. Similar reasoning shows that $b+d$ is always defined and that
the formula for subtraction cannot give an undefined result.
Thus, in interval addition and subtraction         
we achieve the ideal:
{\em Never a Nan}, and this without the need to test.

\subsection{The algorithm for interval multiplication}
If $\ob a,b\cb $ and $\ob c,d\cb $ are bounded, real
intervals, then
$$ \ob a,b\cb  * \ob c,d\cb  = \ob \min(S),\max(S)\cb,$$
where $S = \{a*c, a*d, b*c, b*d\}$.

This formula holds for real rather than floating-point intervals.  It
is several steps away from interval arithmetic.  When we allow the
bounds to be any floating-point number, we introduce the possibility
that they are infinite. In that case we need to be assured that all
four products in $S$ are defined.  Moreover, we want, as much as
possible, to perform only two multiplications, one for each bound.  The
above formula always requires four.

To attain these goals, we classify the intervals $\ob a,b\cb$ and $\ob
c,d\cb$ according to the signs of their elements, as shown in the table
in Figure~\ref{fig:class}.  This classification creates many cases in
which intervals can be multiplied with only one multiplication for each
bound.

\begin{figure}
\begin{tabular}{c||c|c|c|}
Class     & at least one & at least one & Signs of                 \\

of $\ob u,v\cb$   & negative     & positive     & endpoints
        \\
\hline \hline
$M$       &   yes   &   yes     &  $u<0 \wedge v>0$
\\
\hline
$Z$       &   no    &   no      &  $u = 0 \wedge v=0$        \\
\hline
$P$       &   no    &   yes     &  $u \geq 0 \wedge v > 0$        \\
\hline
$P_0$     &   no    &   yes     &  $u = 0 \wedge v > 0$        \\
\hline
$P_1$     &   no    &   yes     &  $u > 0 \wedge v > 0$        \\
\hline
$N$       &   yes   &   no     &    $u < 0 \wedge v \leq 0$      \\
\hline
$N_0$     &   yes   &   no     &    $u < 0 \wedge v = 0$      \\
\hline
$N_1$     &   yes   &   no     &    $u < 0 \wedge v < 0$      \\
\hline
\end{tabular}
\caption{Classification of nonempty intervals according to whether
they contain at least one real of the sign indicated at the top of the
second and third columns. Classes $P$ and $N$ are further decomposed
according to whether they have a zero bound.
As only non-empty intervals are classified, we have $u \leq v$.
}
\label{fig:class}
\end{figure}

The classification yields four cases (for multiplication the
subdivision of $P$ and $N$ do not matter) for each of the operands,
giving at first sight 16 cases. However, when at least one of the
operands classifies as $Z$, several cases collapse. As a result, we are
left with 11 cases.

\begin{theorem}
\label{multTheorem}
If $\ob a,b\cb $ and $\ob c,d\cb $ are real intervals, then
$ \ob a,b\cb  * \ob c,d\cb $
is a real interval whose endpoints are given by the expressions, to
be evaluated as extended reals, in Figure~\ref{multTable1}.
\end{theorem}

\begin{figure*}
\begin{center}
\begin{tabular}{|c|c|c|c|c|}
\hline
Class       & Class      & Left Endpoint     & Right Endpoint & Symmetry\\
of $\ob a,b\cb $  & of $\ob c,d\cb $ & of $\ob
a,b\cb
*\ob c,d\cb $ & of $\ob a,b\cb *\ob c,d\cb $ &
\\
\hline
   P        & P     & $a*c$       & $b*d$          & proved directly \\
   P        & M     & $b*c$       & $b*d$          & proved directly \\
   P        & N     & $b*c$       & $a*d$          & $x*y=-(x*-y)$\\

   M        & P     & $a*d$       & $b*d$          & $x*y=y*x$\\
   M        & M  & $\mn(a*d,b*c)$ & $\mx(a*c,b*d)$ & proved directly \\
   M        & N     & $b*c$       & $a*c$          & $x*y=-(x*-y)$\\

   N        & P     & $a*d$       & $b*c$          & $x*y=-(-x*y)$\\
   N        & M     & $a*d$       & $a*c$          & $x*y=-(-x*y)$\\
   N        & N     & $b*d$       & $a*c$          & $x*y=-(x*-y)$\\

   Z        & P,M,N,Z &   0         &   0          & proved directly\\
  P,M,N       & Z     &   0         &   0          & proved directly \\
\hline
\end{tabular}
\caption{Case analysis for multiplication of real intervals,
$\ob a,b\cb *\ob c,d\cb $.}
Results for floating-point intervals are obtained by performing the
lower-bound (upper-bound) computations rounded toward $-\infty$
($+\infty$).
\label{multTable1}
\end{center}
\end{figure*}

In  \cite{hckvnmdn01} the cases indicated as such in the table in
Figure~\ref{multTable1} are proved directly. The other cases can be
proved by symmetry from the case proved already. The symmetries
applied are based on the identities $x*y=-(x*-y)$ or 
similar ones shown in the last column in the table.

The proofs first show the correctness of the scalar products for
bounded real intervals. To allow for floating-point intervals, which
can be unbounded, we have to consider whether the products are
defined.  Let us consider as example the top line according to which
$\ob a,b \cb * \ob c,d \cb = \ob a*c, b*d \cb$.  The undefined cases
occur when one operand is 0 and the other $\infty$. It is possible for
$a$ or $c$ to equal 0, but neither can be infinite: because of the
classification $P$, they cannot be $-\infty$; because of their being
lower bounds, they cannot be $+\infty$.

Let us now consider $b*d$. It is possible for $b$ or $d$ to equal
$+\infty$, but neither can be 0 because of the classification $P$.  One
may verify that in every case of the table in Figure~\ref{multTable1}
undefined values are avoided by a combination of definitions
\ref{brackets} and \ref{flptOrReal} and the classification of the case
concerned.

We need tests to identify the right case in the table anyway to
minimize the number of multiplications. We obtain as a bonus the
saving of tests to avoid undefined values.
Thus, in interval multiplication
we achieve the ideal:
{\em Never a Nan}, and this without the need to test.

\subsection{Division}

For interval multiplication the classification of the interval
operands in the classes $P$, $M$, $N$, and $Z$ is sufficient. For
interval division it turns out that the further subdivision of $P$
into $P_0$ and $P_1$ and of $N$
into $N_0$ and $N_1$ (see the table in figure~\ref{fig:class})
is relevant for the dividend.

\begin{theorem}
\label{divTheorem}
If $\ob a,b\cb $ and $\ob c,d\cb $ are real intervals, then $ \ob
a,b\cb  \oslash \ob c,d\cb $ is the least floating-point interval
containing the real interval whose endpoints are given as the ``general
formula'' column in Figure~\ref{divTable1} unless the specified
condition in the next column holds, in which case the result is given
by the exception case in column 5.

\end{theorem}

\begin{figure*}
\begin{center}
\begin{tabular}{|c|c|l|c|l|c|}
\hline
Class    & Class    & $\ob a,b\cb \oslash \ob c,d\cb $     &        & $\ob a,b\cb /\ob c,d\cb $  &  \\
of $\ob a,b\cb $ & of $\ob c,d\cb $ & general formula   &unless  & exception case &  \\
\hline
  $P_1$     & $P$   & $\ob a/d,b/c\cb  \setminus \{0\}$  & $c=0$   & $\ob a/d,\infty\cb \setminus\{0\}  $ &$D$ \\
  $P_0$     & $P$   & $\ob 0,b/c\cb $                    & $c=0$   & $\ob -\infty,\infty\cb                   $ &$D$ \\
  $M$       & $P$   & $\ob a/c,b/c\cb $                  & $c=0$   & $\ob -\infty,\infty\cb             $ &$D$ \\
  $N_0$     & $P$   & $\ob a/c,0\cb $                    & $c=0$   & $\ob -\infty,\infty\cb                  $ &$S_2$\\
  $N_1$     & $P$   & $\ob a/c,b/d\cb  \setminus \{0\}$  & $c=0$   & $\ob -\infty,b/d\cb \setminus \{0\}$ &$S_2$\\

  $P_1$     & $M$   & $(\ob -\infty,a/c\cb \cup\ob a/d,\infty\cb )\setminus \{0\}$  &         &           &$D$ \\
  $P_0$     & $M$   & $\ob -\infty,+\infty\cb $          &         &                                &$D$ \\
  $M$       & $M$   & $\ob -\infty,+\infty\cb $          &         &                                &$D$ \\
  $N_0$     & $M$   & $\ob -\infty,+\infty\cb $          &         &                                &$S_2$\\
  $N_1$     & $M$   & $(\ob -\infty,b/d\cb \cup\ob b/c,\infty\cb ) \setminus \{0\}$  &         &                                &$S_2$\\

  $P_1$     & $N$   & $\ob b/d,a/c\cb  \setminus \{0\}$  & $d=0$   & $\ob -\infty,a/c\cb \setminus\{0\} $ &$S_1$\\
  $P_0$     & $N$   & $\ob b/d,0\cb $                    & $d=0$   & $\ob -\infty,\infty\cb                  $ &$S_1$\\
  $M$       & $N$   & $\ob b/d,a/d\cb $                  & $d=0$   & $\ob -\infty,\infty\cb             $ &$S_1$\\
  $N_0$     & $N$   & $\ob 0,a/d\cb $                    & $d=0$   & $\ob -\infty,\infty\cb                   $ &$S_2$\\
  $N_1$     & $N$   & $\ob b/c,a/d\cb  \setminus \{0\}$  & $d=0$   & $\ob b/c,\infty\cb \setminus\{0\}  $ &$S_2$\\
  $Z$     & $P_1,N_1$ & $\ob 0,0 \cb$                       &      & 
                            &     \\
  $Z$ & $P_0,M,N_0,Z$ & $\ob -\infty,+\infty \cb$           &      & 
                            &     \\
  $P_1,N_1$ & $Z$     & $\emptyset$                         &      &
                            &     \\
  $P_0,M,N_0,Z$ & $Z$ & $\ob -\infty,+\infty \cb$           &      &
                            &     \\
\hline
\end{tabular}
\caption{Case analysis for relational division of real intervals, $\ob
a,b\cb /\ob c,d\cb $ when $a \leq b$, $c \leq d$.
The last column refers to how the formula has been
proved (``$D$'' for a direct proof, ``$S_1$'' and ``$S_2$'' refer to a
symmetry used to reduce it to an earlier case.) The ``class'' labels,
$N,N_1,N_0,M,P_0,P_1,P$ are as in Figure \ref{fig:class}.
}
\label{divTable1}
\end{center}
\end{figure*}

In  \cite{hckvnmdn01} the cases indicated as such in the table in
Figure~\ref{divTable1} are proved directly. The other cases can be
proved by symmetry from the case proved already. The symmetries
used are based on the identities $x/y = -(x/-y)$ (indicated as $S_1$)
and $x/y = -(-x/y)$ (indicated by $S_2$).

The proofs first show the correctness of the scalar products for
bounded real intervals. To allow for floating-point intervals, which
can be unbounded, we have to consider whether the products are defined.
In the column labelled ``unless'' we find the values for which an
undefined value occurs. In the ``exception case'' column we find the
correct value for the exception case.
{\em In every case, evaluating the formula in the third column in IEEE
standard floating-point arithmetic in the exception case is defined
and gives the infinity of the right sign,} as shown in column 5.
This property depends on a zero
lower bound being $+0$ and a zero upper bound being $-0$, as required by
definition~\ref{flptOrReal}.

Let us now consider potentially undefined cases. In case of division
these are $\infty/\infty$ and $0/0$. Consider for example the top 
line according to which
$\ob a,b \cb \oslash \ob c,d \cb = \ob a/d, b/c \cb \setminus \{0\}$.
Because of the classification $P_1$, $a$ can be neither infinite nor
zero. This ensures that $a/d$ is defined.
Because of the $P_1$ classification, $b$ cannot be zero.
It is possible for $b$ to be infinite, but not for $c$ because of the
$P$ classification.
This ensures that $b/c$ is defined.

One may verify that in every case of the table in
Figure~\ref{divTable1} undefined values are avoided by a combination
of definition~\ref{flptOrReal} and the classification of the case concerned.
Thus, in relational interval division
we achieve the ideal:
{\em Never a Nan}, and this without the need to test.

\section{Related work}
For most of the time since the beginning of interval arithmetic, two
systems have coexisted.
One was the official one, where intervals were bounded, and division
by an interval containing zero was undefined.
Recognizing the unpracticality of this approach, there was also a
definition of ``extended'' interval arithmetic
\cite{khn68}
where these limitations
were lifted. Representative of this state of affairs are the
monographs by Hansen \cite{hnsn92} and Kearfott \cite{krftt96}.
However, here the specification of interval division is quite far from an
efficient implementation that takes advantage of the IEEE floating-point
standard. The specification is indirect via multiplication by
the interval inverse. There is no consideration of the possibility of
undefined operations: presumably one is to perform a test before each
operation.

Steps beyond this were taken by Older \cite{ldr89ias} in connection
with the development of BNR Prolog. A different approach has been taken
by Walster \cite{walster98}, who pioneered the idea that intervals are
sets of values rather than abstract elements of an interval algebra.
In Walster shares our objective to obtain a closed system of arithmetic
without exceptions. He attains this objective in a different way: by
including the infinities among the possible values of the variables. In
our approach, the variables can only take reals as values; the
infinities are only used for the representation of unbounded sets of
reals. In this way, the conventional framework of calculus, where
variables are restricted to the reals, needs no modification.

\section{Conclusions}

We have presented the result of some recent developments in interval
arithmetic that lead to a system with the following properties.

\begin{itemize}
\item
{\em Correctness}
The interval operations are such that their result includes
all real numbers that are
possible as values of the variables according to the mathematical
model.

\item
{\em Freedom of exceptions}
No floating-point operation needs raise an exception.
All divisions by zero are
defined and give the correct result: an infinity of the correct
sign. This is achieved by a  zero lower (upper) bound being $+0$
($-0$). Mathematically speaking, the system is a closed interval
algebra. We do not emphasize the algebra aspect, because it is not
important whether it has any interesting properties.
Other approaches have limited the applicability of interval arithmetic
in their pursuit of a presentable algebra.

\item
{\em Efficiency}
The system is efficient in that tests are only needed to determine the
right case in the tables in Figures \ref{multTable1} and
\ref{divTable1}.  Tests are not necessary to avoid exceptions.

\end{itemize}

\noindent
These properties lead to several observations about the
floating-point standard from the point of view of interval arithmetic:

\begin{itemize}
\item
{\em Exceptions}
Freedom from exceptions has interesting implications for the standard.
A considerable part of the definition effort, and presumably also of
the implementation effort, is concerned with defining, signaling, and
trapping exceptions caused by overflow, underflow and undefined
operations.  A processor where the floating-point arithmetic is
interval arithmetic can omit this as unnecessary ballast.

Let us review the five exceptions. {\em Invalid Operation} is
prevented by the design of the algorithms.
{\em Division by Zero} does occur in our interval arithmetic and is
designed to yield the correct result. So it should not be an
exception.
{\em Overflow} occurs in the sense that a real $x$ can result in real
arithmetic such that $\phi(x)$ is the interval between the greatest
finite floating-point number and $+\infty$. This result is
mathematically correct and therefore the desired one.
There is no reason to terminate computation: it should
not be an exception.
{\em Underflow} means that a lower bound zero is substituted for a
nonzero bound with very small absolute value. This is correct and no
reason to terminate computation.
{\em Inexact} result: this might be of some use, but is certainly not
essential for interval arithmetic.

\item
{\em Signed zeros}
Often signed zeros are regarded as an unavoidable, but regrettable
artifact of the sign-magnitude format of floating-point numbers. 
It is fortunate that the drafters of the standard have nonetheless
taken them seriously and defined sensible conventions for operations
involving zeros. Especially having the right sign of a zero bound
turns out to be useful in interval division.

\item
{\em Denormalized numbers}
of view of interval arithmetic, denormalized numbers seem to be neither
useful nor harmful.  It is different from the point of view of interval
constraints.  This a method \cite{bnldr97,vhlmyd97} for using interval
arithmetic to solve systems of constraints with real-valued variables.
Interval arithmetic is used for the basic operations in constraint
propagation. This is an iteration that can be slowed down by
denormalized numbers when the limit is zero, even when operations on
denormalized numbers are performed at normal speed. Thus the presence
of denormalized numbers only plays a role as a performance bug that occurs
gratuitously, and fortunately rarely, in this special case.

An argument that is advanced in favour of denormalized numbers is that
it justifies compiler optimizations that rely on
certain mathematical equivalences that hold only in the presence of
denormalized numbers.
This is of no interest from the point of view of interval constraints.
Any mathematically correct transformation can be performed on
the set of constraints without changing the set of solutions obtained
by a correctly implemented interval constraint system.
This correctness is not dependent on the presence of denormalized
numbers. In fact, it only depends on the finite floating-point numbers
being {\em some} subset $F$ of the reals, as described in this paper. 
Because of this independence, elaborate symbolic processing far beyond
currently contemplated compiler optimizations is taken for granted in
interval constraints.
\end{itemize}

\section{Acknowledgments}
Many thanks to Belaid Moa for pointing out errors.
We acknowledge generous support from the Natural Science and
Engineering Research Council NSERC.

\end{document}